\documentstyle[12pt]{article}
\begin{document}
\title{THE ORIGIN OF A METRIC}
\author{B.G.Sidharth \\
B.M.Birla Science Centre,Hyderabad, 500463,India}
\date{}
\maketitle
\begin{abstract}
In the context of earlier work, we investigate the emergence of a "distance"
in the physical world.
For this we consider a Cantor ternary like process, but much more general:
properties like perfectness and disconnectedness are not invoked, but instead
we deal with Borel sets. An interesting case from a physical point of view is
considered: when the process is truncated.
\end{abstract}
\section{The Origin of a Metric}
We first makes a few preliminary remarks. When we talk of a metric
or the distance between two "points" or "particles", a concept
that is implicit is that of topological "nearness" - we require an
underpinning of a suitably large number of "open" sets\cite{r1}.
Let us now abandon the absolute or background space time and
consider, for simplicity, a universe (or set) that consists solely
of two particles. The question of the distance between these
particles (quite apart from the question of the observer) becomes
meaningless. Indeed, this is so for a universe consisting of a
finite number of particles. For, we could isolate any two of them,
and the distance between them would have no meaning. We can
intuitiively appreciate that we would infact need distances of
intermediate or more generally, other points.\\ In earlier
work\cite{r2}, motivated by physical considerations we had
considered a series of nested sets or neighbourhoods which were
countable and also whose union was a complete Hausdorff space. The
Urysohn Theorem was then invoked and it was shown that the space
of the subsets was metrizable. The argument went something like
this.\\ In the light of the above remarks, the concepts of open
sets, connectedness and the like reenter in which case such an
isolation of two points would not be possible.\\ More formally let
us define a neighbourhood of a particle (or point or element) $A$
of a set of particles as a subset which contains $A$ and atleast
one other distinct element. Now, given two particles (or points)
$A$ and $B$, let us consider a neighbourhood containing both of
them, $n(A,B)$ say. We require a non empty set containing atleast
one of $A$ and $B$ and atleast one other particle $C$, such that
$n(A,B) \subset n(A,C)$, and so on. Strictly, this "nested"
sequence should not terminate. For, if it does, then we end up
with a set $n(A,P)$ consisting of two isolated "particles" or
points, and the "distance" $d(A,P)$ is meaningless.\\ We now
assume the following property\cite{r2}: Given two distinct
elements (or even subsets) $A$ and $B$, there is a neighbourhood
$N_{A1}$ such that $A$ belongs to $N_{A1}$, $B$ does not belong to
$N_{A1}$ and also given any $N_{A1}$, there exists a neighbourhood
$N_{A_\frac{1}{2}}$ such that $A \subset N_{A_\frac{1}{2}} \subset
N_{A1}$, that is there exists an infinite topological closeness.\\
>From here, as in the derivation of Urysohn's lemma\cite{r1}, we
could define a mapping $f$ such that $f(A) = 0$ and $f(B) = 1$ and
which takes on all intermediate values. We could now define a
metric, $d(A,B) = |f(A) - f(B)|$. We could easily verify that this
satisfies the properties of a metric.\\ With the same motivation
we will now deduce a similar result, but with different
conditions. In the sequel, by a subset we will mean a proper
subset, which is also non null, unless specifically mentioned to
be so. We will also consider Borel sets, that is the set itself
(and its subsets) has a countable covering with subsets. We then
follow a pattern similar to that of a Cantor ternary set
\cite{r1,r3}. So starting with the set $N$ we consider a subset
$N_1$ which is one of the members of the covering of $N$ and
iterate this process so that $N_{12}$ denotes a subset belonging
to the covering of $N_1$ and so on.\\ We note that each element of
$N$ would be contained in one of the series of subsets of a sub
cover. For, if we consider the case where the element $p$ belongs
to some $N_{12\cdots m}$ but not to $N_{1,2,3\cdots m+1}$, this
would be impossible because the latter form a cover of the former.
In any case as in the derivation of the Cantor set, we can put the
above countable series of sub sets of sub covers in a one to one
correspondence with suitable sub intervals of a real interval
$(a,b)$.\\ {\large \bf{Case I}}\\ If $N_{1,2,3\cdots m} \to$ an
element of the set $N$ as $m \to \infty$, that is if the set is
closed, we would be establishing a one to one relationship with
points on the interval $(a,b)$ and hence could use the metric of
this latter interval, as seen earlier.\\ {\large \bf{Case II}}\\
It is interesting to consider the case where in the above
iterative countable process, the limit does not tend to an element
of the set $N$, that is set $N$ is not closed and has what we may
call singular points. We could still truncate the process at
$N_{1,2,3\cdots m}$ for some $m > R$ arbitrary and establish a one
to one relationship between such truncated subsets and arbitrarily
small intervals in $a,b$. We could still speak of a metric or
distance between two such arbiitrarily small intervals.\\ This
case is of interest because of recent work which describes
elementary particles as, what may be called Quantum Mechanical
Kerr-Newman Black Holes or vortices, where we have a length of the
order of the Compton wavelength (that is $10^{-12}cms$ or less),
within which spacetime as we know it breaksdown. Such cut offs
lead to a non commutative geometry and what may be called fuzzy
spaces\cite{r4},\cite{r5},\cite{r6},\cite{r7},
\cite{r8},\cite{r9}.(We note that the centre of the vortex is a
singular point). In any case, the number of particles in the
universe is of the order $10^{80}$, which approxiimates infinity
from a physicist's point of view.\\ {\large \bf{Remarks}}\\
Interestingly, we usually consider two types of infinite sets -
those with cardinal number $n$ corresponding to countable
infinities, and those with cardinal number $c$ corresponding to a
continuum, there being nothing inbetween. This is the well known
but unproven Continuum Hypotheses.\\ What we have shown with the
above process is that it is possible to concieve an intermediate
possibility with a cardinal number $n^p, p > 1$.\\ We also note
the similarity with transfinite Cantor sets. But in this latter
case three properties are important: the set must be closed i.e.
it must contain all its limit points, perfect i.e. in addition
each of its points must be a limit point and disconnected i.e. it
contains no nonnull open intervals. Only the first was invoked in
Case I.\\ Finally we will remark on an origin for spatial
dimensions, in the context of the above elementary particle model
and an associated consistent cosmology\cite{r10}, in which given
$n$ particles, $\sqrt{n}$ would be fluctuationally created from a
background Quantum vaccuum within the undefined time interval at
the Compton scale referred to above. If there are $n$ points,
these could be lined up as a single dimension. But if with each
point other points are associated, not belonging to the origiinal
set then this would define another dimension. That is the case
with the fluctuationally created particles: There would now be not
$n$ but $n^{3/2}$ particles or points. Further with each of the
fluctuationally created particles, there would be $n^{1/4}$
further fluctuationally created particles\cite{r11}. So in effect
we would have to deal with not $n^{3/2}$ but $n^{7/4}$ points and
so on so that finally with the $n$th point or particle we would
have $n^2$ points. The total number of points would therefore be
$\Sigma n^2$ or $\sim n^3$ giving us the three spatial dimensions.

\end{document}